\documentclass{amsart}
\usepackage{caption,hyperref}
\usepackage{enumerate,amsmath,amssymb,graphicx}
\usepackage{subfig,color}

\def\NN{\mathbb N}

\newtheorem{theorem}{Theorem}

\newtheorem{remark}{Remark}

\title{Solutions of Bernoulli equations in the fractional setting}
\author{Mirko D'Ovidio, Anna Chiara Lai,  Paola Loreti}
\address{ Sapienza Università di Roma, Dipartimento di Scienze di Base e Applicate per l'Ingegneria; }
\email{\{mirko.dovidio,annachiara.lai,paola.loreti\}@uniroma1.it}

\begin{document}
\begin{abstract}{
We present a general  series representation formula for the local solution of  Bernoulli equation with Caputo fractional derivatives.
  We then focus on a generalization of the fractional logistic equation and we present some related numerical simulations. }
\end{abstract}
\maketitle

\section{Introduction}
Interest for time fractional  evolutive  systems  has progressively grown in recent years: models   arising  in nature with aspects  related to non-local behaviour need a study in the fractional setting, see \cite{pod, fracmod} for an overview and  \cite{fracmod2, fracmod} for fractional growth models for social and biological dynamics. The fractional derivatives are indeed non-local operators, that is convolution-type operators. In the applied sciences, the main interest in fractional models is due to the fact that such models introduce the so-called memory effect. This effect is mainly justified by the non-locality of the time-fractional derivative and it seems to be relevant in the characterization of many applied models. A second reading is given in terms of the delaying effect. Indeed, the time-fractional derivative introduces a different clock for the underlying model as in case of the relaxation equation. When the order of the fractional derivative is 1 the underlying model emerges.

Here we locally solve the following Cauchy system, involving a fractional Bernoulli equation of the form
\begin{equation}D_t^\beta u+ a_0 u = a_1u^{p+1}, \qquad u(0) = u_0\label{eq1}  \end{equation}
where $p\in\NN$,  with $p\geq 1$, $a_0$ and $a_1$ real numbers and $D^\beta_t$ denotes the Caputo derivative.
If $\beta=1$ then $D^\beta_t u= u'$
and \eqref{eq1} is the Bernoulli equation, studied by Jacob Bernoulli (1695). We remark that Bernoulli equations (with $\beta=1$) arise in non-linear models of production and capital accumulation, in particular when polynomial production functions are considered, see \cite[Chapter 6.3]{H}. As a particular case, the exact solution in the case $a_0=a_1=-1$ and $\beta=1$ is given by
$$u(t)=\frac{e^t} {\sqrt[p]{c_0+e^{pt}}}; \qquad \text{where } c_0:= u_0^{-p}-1.$$

\begin{figure}
\centering \includegraphics[scale=.45]{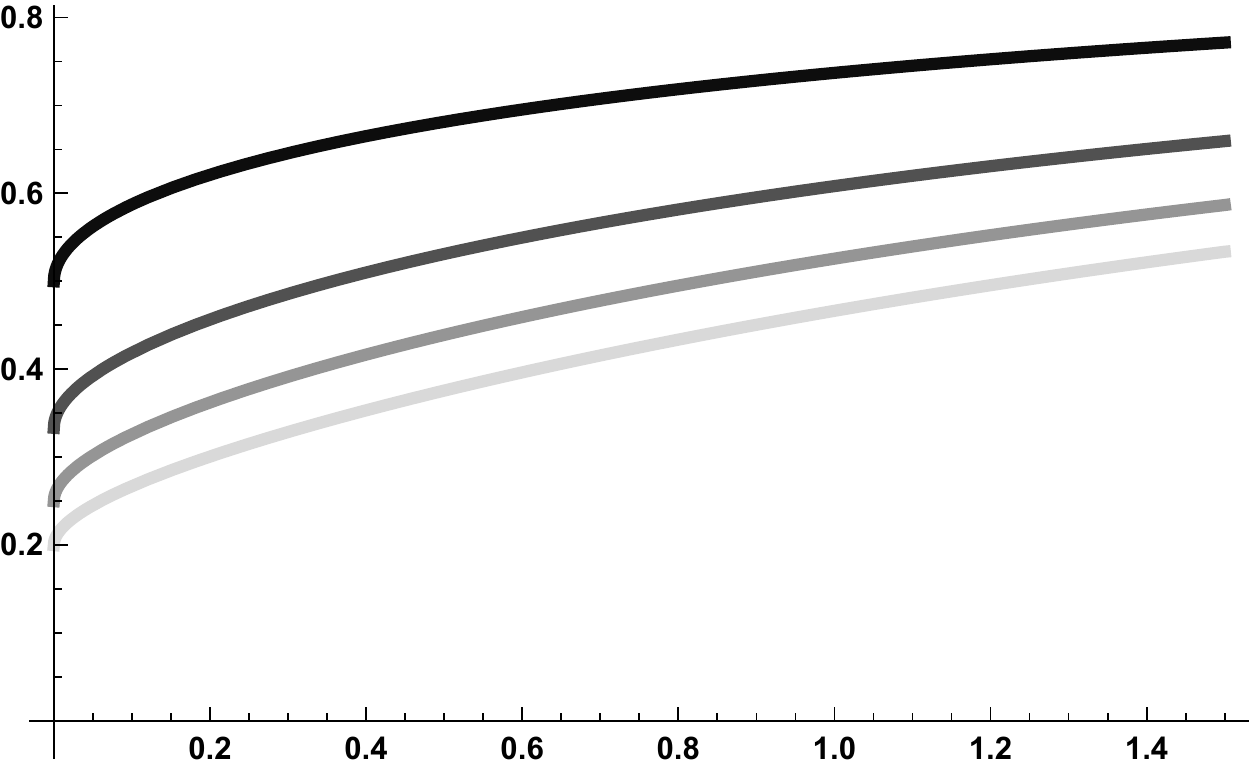}\quad\includegraphics[scale=.45]{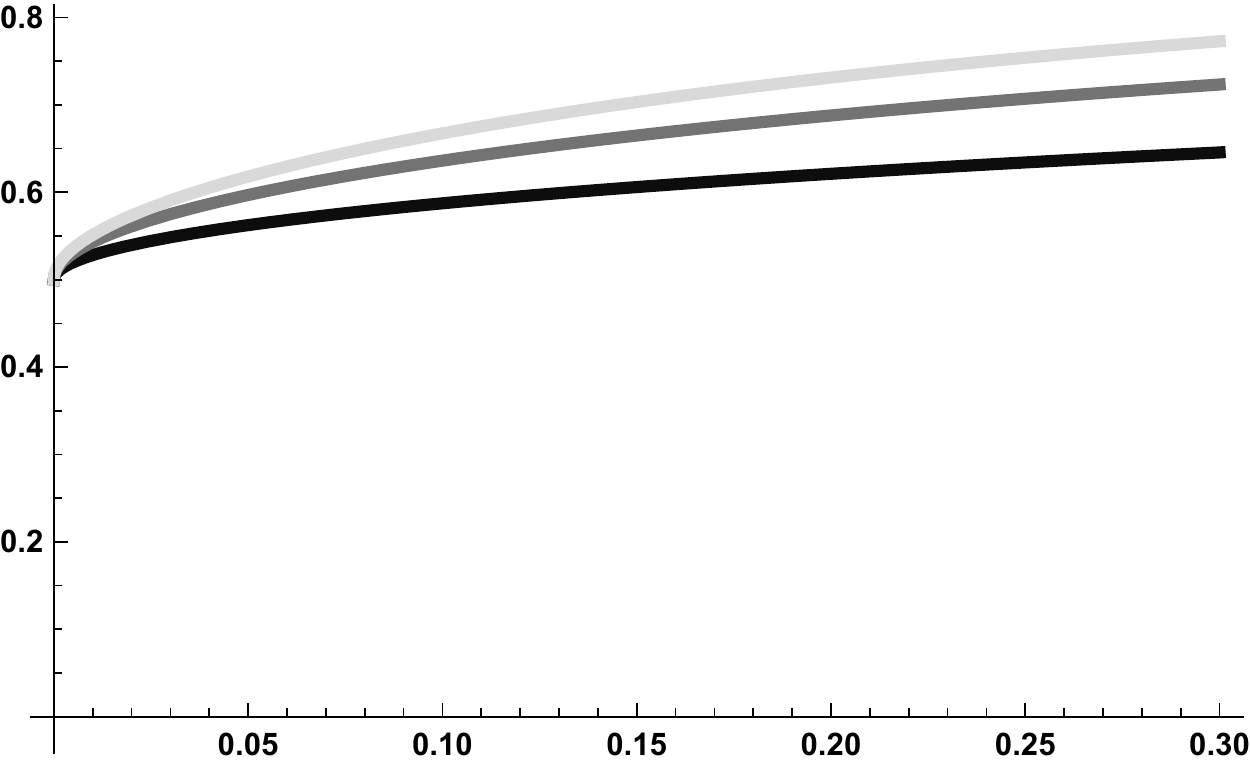}\
\caption{On the left, solutions of the logistic equation $D^\beta u-u=-u^2$ (on the left)  with $\beta=1/2$ and $u_0=1/2,1/3,1/4,1/5$; darker lines correspond to bigger $u_0$s. On the right, solutions of the logistic equation $D^\beta u-u=-u^{p+1}$ with $\beta=u_0=1/2$ and $p=1,2,3$; darker lines correspond to bigger $p$s.\label{f1}
}
\end{figure}

Going to the fractional setting, we have that similar approaches cannot be followed.
As it is well known, also the solution of the fractional logistic equation --corresponding to $p=1$ and $a_0=a_1=-1$ in \eqref{eq1}-- was an open problem and in \cite{DL} the first and the third author were
able to solve the fractional logistic equation by series representation, giving a detailed formula involving Euler
numbers for $u_0 =1/2$. This approach was then applied to SIS epidemic models in \cite{BDL} and also further investigated in \cite{nieto}. The present study extends the result in \cite{DL} to general initial data and to Bernoulli equations of general degree $p+1$: we present a recursive formula for the coefficients of the solutions and explicit closed formulas for the first terms. Note that  the relation with Euler numbers for general initial data, even in the logistic case $p=1$, appears to be lost, but  the general recursive formula preserves its structure, based on generalized binomial coefficients that were introduced in \cite{DL} and further  investigated in \cite{DLL}. Then the proposed method is applied to the particular case $a_0=a_1=\pm1$, related to the fractional logistic equation and we present a qualitative analysis of the solutions based on numerical simulations -- see Figure \ref{f1}.

\subsection{Prelimininaries on fractional calculus}

Let us consider the set $AC([a,b])$ of continuous functions with derivative in $L_1([a,b])$. Thus, $v\in AC([a,b])$ is continuous and such that $v^\prime = g \in L_1([a,b])$, that is $v$ has the representation
\begin{align}
v(t)= v(a) + \int_a^t g(s)ds, \quad t \in [a,b].
\end{align} 
We notice that the space $AC([a,b])$ coincides with the Sobolev space
\begin{align*}
W^{1,1}([a,b]) = \{v \in L_1([a,b])\, :\, v^\prime \in L_1([a,b])\}
\end{align*}
endowed with the norm
\begin{align*}
\|v\|_{W^{1,1}} = \|v\|_{L_1} + \|v^\prime \|_{L_1}.
\end{align*}
For $v\in AC([a,b])$ and $\beta \in (0,1)$ we introduce the Riemann-Liouville derivative of $v$,
\begin{align}
\label{DefRL}
\mathcal{D}^\beta_t v(t) : = \frac{1}{\Gamma(1-\beta)} \frac{d}{dt} \int_a^t v(s) (t-s)^{-\beta} ds, \quad t \in [a,b]
\end{align}
and the Caputo-Djarbashian derivative of $v$,  
\begin{align}
\label{DefCJ}
D^\beta_t v(t) : = \frac{1}{\Gamma(1-\beta)} \int_a^t v^\prime(s) (t-s)^{-\beta} ds, \quad t \in [a,b].
\end{align}
Further on we use the following relation between derivatives
\begin{align}
\mathcal{D}^\beta_t v(t) = & \frac{1}{\Gamma(1-\beta)}  \frac{v(a)}{(t-a)^\beta} + D^\beta_t v(t). \label{RelationRLC}
\end{align}
The relation \eqref{RelationRLC}, together with the existence of the  derivatives \eqref{DefRL} and \eqref{DefCJ}, hold a. e. on $[a,b]$ and $\mathcal{D}^\beta_t v \in L_q([a,b])$ with $1\leq q \leq 1/\beta$ (see for example \cite[page 28]{Diethelm}). 

%
%
%

We consider throughout fractional equations on $[0, b)$. Let us underline that, if $v(0)=0$, then formula \eqref{RelationRLC} gives the equivalence
\begin{align*}
\mathcal{D}^\beta_t v = D^\beta_t v.
\end{align*} 

\section{Fractional Bernoulli equations}

{

Let us introduce 
\begin{align}
\label{solEQ}
u(t) = \sum_{n \geq 0} c_n^{(1)}\, \frac{t^{\beta n}}{\Gamma(\beta n + 1)}, \quad t \in (0, r)
\end{align}
where  
\begin{align*}
\left[\!\begin{matrix} \, n \, \\ \, k\,  \end{matrix}\!\right]_\beta = \frac{\Gamma(n\beta +1)}{\Gamma(k\beta +1) \, \Gamma((n-k)\beta + 1)}, \quad 0\leq k \leq n, \quad k,n \in \mathbb{N}, \quad \beta \in (0,1).
\end{align*}
is a generalized binomial coefficient, $r>0$ is the radius of convergence and $c_n^{(1)}$ are real coefficients.

\begin{theorem}
The unique continuous solution on $[0,b) \subset [0, r)$ to
\begin{align*}
D^\beta_t u + a_0\, u = a_1\, u^{p+1}, \quad u(0)=u_0, \quad p \in \mathbb{N}, \quad a_0,a_1 \in \mathbb R, \quad \beta \in (0,1)
\end{align*}
has the series representation \eqref{solEQ} on $[0,b)$ where
\begin{align}
& c_0^{(1)} = u_0\notag\\
& c_{n+1}^{(1)} = -a_0\, c_{n}^{(1)} + a_1 \sum_{k=0}^n \left[\!\begin{matrix} \, n \, \\ \, k\,  \end{matrix}\!\right]_\beta\, c^{(p)}_k\, \, c^{(1)}_{n-k} \qquad \text{for } n \geq 0\label{c1n}
\end{align}
\end{theorem}

}
\begin{proof}

We consider 
\begin{align*}
D^\beta_t u +  {a_0}u = {a_1}u^{p+1}, \quad u(0)=u_0
\end{align*}
where $D^\beta_t$ denotes the Caputo derivative. 
Let us denote the power of $u$ as follows
\begin{align}
u(t) = \sum_{n \geq 0} c^{(1)}_n \frac{t^{\beta n}}{\Gamma(\beta n +1)}\label{series}\\
\notag\\
u^2(t) = \sum_{n\geq 0} c^{(2)}_n\, \frac{t^{\beta n}}{\Gamma(\beta n +1)}\notag
\end{align}
\begin{align*}
u^3(t) = \sum_{n\geq 0} \left( \sum_{k=0}^n \left[\!\begin{matrix} \, n \, \\ \, k\,  \end{matrix}\!\right]_\beta\, c^{(2)}_k\, \, c^{(1)}_{n-k} \right) \frac{t^{\beta n}}{\Gamma(\beta n +1)} = \sum_{n \geq 0} c^{(3)}_n \, \frac{t^{\beta n}}{\Gamma(\beta n +1)}
\end{align*}
where
\begin{align*}
c^{(2)}_k = \sum_{s=0}^k \left[\!\begin{matrix} \, k \, \\ \, s\,  \end{matrix}\!\right]_\beta\, c^{(1)}_s\, \, c^{(1)}_{k-s} \quad \Rightarrow \quad c^{(3)}_n = \sum_{k=0}^n \left[\!\begin{matrix} \, n \, \\ \, k\,  \end{matrix}\!\right]_\beta\, c^{(2)}_k\, \, c^{(1)}_{n-k}
\end{align*}
and, by further iterations,
\begin{align*}
u^{p+1}(t) = \sum_{n \geq 0} c^{(p+1)}_n \frac{t^{\beta n}}{\Gamma(\beta n +1)}\\
\end{align*}
where
\begin{align}\label{cnp}
c^{(p+1)}_n = \sum_{k=0}^n \left[\!\begin{matrix} \, n \, \\ \, k\,  \end{matrix}\!\right]_\beta\, c^{(p)}_k\, \, c^{(1)}_{n-k}.
\end{align}
{Note that for $p\in\mathbb N$
$$c_0^{(p)}=\left(c_0^{(1)}\right)^p.$$
}
Therefore
\begin{align*}
{a_0}u - {a_1}u^{p+1} 
= & u({a_0}-{a_1}u^p) = \sum_{i} c^{(1)}_i \frac{t^{\beta i}}{\Gamma(\beta i + 1)} \cdot \left( {a_0} - {a_1}\sum_j c_j^{(p)} \frac{t^{\beta j}}{\Gamma(\beta j + 1)} \right)\\
= &{a_0} \sum_{i} c^{(1)}_i \frac{t^{\beta i}}{\Gamma(\beta i + 1)} - {a_1}\sum_{i} c^{(1)}_i \frac{t^{\beta i}}{\Gamma(\beta i + 1)} \sum_j c_j^{(p)} \frac{t^{\beta j}}{\Gamma(\beta j + 1)}\\
= & {a_0}\sum_j c_j^{(1)} \frac{t^{\beta j}}{\Gamma(\beta j + 1)}  - {a_1}\sum_n c_n^{(p+1)} \frac{t^{\beta n}}{\Gamma(\beta n + 1)}\\
= & \sum_n \big( {a_0}c_n^{(1)} -{a_1} c_n^{(p+1)} \big)\frac{t^{\beta n}}{\Gamma(\beta n + 1)}\\
= & \sum_n \left({a_0} c_n^{(1)} -{a_1} \sum_{k=0}^n \left[\!\begin{matrix} \, n \, \\ \, k\,  \end{matrix}\!\right]_\beta\, c^{(p)}_k\, \, c^{(1)}_{n-k} \right)\frac{t^{\beta n}}{\Gamma(\beta n + 1)}.
\end{align*}
On the other hand, from \eqref{RelationRLC}, we have that
\begin{align*}
D^\beta_t u(t) = - u_0\, \frac{t^{-\beta}}{\Gamma(1-\beta)} + \mathcal{D}^\beta_t \sum_{n \geq 0} c_n^{(1)} \frac{t^{\beta n}}{\Gamma(\beta n +1)}
\end{align*}
where, after some calculation, from \eqref{DefRL},
\begin{align*}
\mathcal{D}^\beta_t \frac{t^{\beta n}}{\Gamma(\beta n +1)} = \frac{t^{\beta(n-1)}}{\Gamma(\beta (n-1)+1)}, \quad t \in [0, b).
\end{align*}
Thus, we obtain
\begin{align*}
D^\beta_t u(t) 
= & - u_0\, \frac{t^{-\beta}}{\Gamma(1-\beta)} + \sum_{n\geq 0} c_n^{(1)}\, \frac{t^{\beta(n-1)}}{\Gamma(\beta (n-1)+1)}\\
= & - u_0\, \frac{t^{-\beta}}{\Gamma(1-\beta)} + c_0^{(1)} \frac{t^{-\beta}}{\Gamma(-\beta +1)} + \sum_{n\geq 1} c_n^{(1)}\, \frac{t^{\beta(n-1)}}{\Gamma(\beta (n-1)+1)}.
\end{align*}
From the fact that $u_0=c_0^{(1)}$ by construction, we write
\begin{align*}
D^\beta_t u(t) = \sum_{n \geq 0} c_{n+1}^{(1)} \frac{t^{\beta n}}{\Gamma(\beta n + 1)}, \quad t \in [0, b].
\end{align*}
Then the solution to
\begin{align*}
D^\beta_t u = - {a_0}u + {a_1}u^{p+1} = - u({a_0}-{a_1}u^p)
\end{align*}
can be written in terms of the coefficients $c^{(1)}_n$, $n \in \mathbb{N}$ given by
\begin{align*}
c^{(1)}_{n+1}  &=-{a_0 c_n^{(1)}+ {a_1}c_{n}^{(p+1)}}\\
 &{= -c_n^{(1)} \big({a_0} -{a_1} c^{(p)}_0 \big) +{a_1}  \sum_{k=1}^n \left[\!\begin{matrix} \, n \, \\ \, k\,  \end{matrix}\!\right]_\beta\, c^{(p)}_k\, \, c^{(1)}_{n-k}},\notag\\
&=-{a_0} c_n^{(1)}+ {a_1} \sum_{k=0}^n \left[\!\begin{matrix} \, n \, \\ \, k\,  \end{matrix}\!\right]_\beta\, c^{(p)}_k\, \, c^{(1)}_{n-k} \qquad \text{for } n \geq 0.
\end{align*}
\end{proof}
\subsection{Some closed formulas}
{
The first few element of the sequence $c^{(p)}_n$, $p \in \mathbb{N}$ (over the index $p$) are
\begin{align*}
& c^{(1)}_n, \quad n \in \mathbb{N}\\
\\
& c^{(2)}_n = \sum_{s=0}^n \left[\!\begin{matrix} \, n \, \\ \, s\,  \end{matrix}\!\right]_\beta\, c^{(1)}_s\, \, c^{(1)}_{n-s}\\
\\
& c^{(3)}_n = \sum_{k=0}^n \left[\!\begin{matrix} \, n \, \\ \, k\,  \end{matrix}\!\right]_\beta\, \sum_{s=0}^k \left[\!\begin{matrix} \, k \, \\ \, s\,  \end{matrix}\!\right]_\beta\, c^{(1)}_s\, \, c^{(1)}_{k-s}\, c^{(1)}_{n-k}.
\end{align*}

We compute the first terms of $c^{(1)}_n$. Fix $p\geq 1$ and note that $c^{(1)}_0= u_0$ and if $n=1$ then 
\begin{equation}c^{(1)}_1= -u_0({a_0}-{a_1}u_0^p).\label{c11}\end{equation}

To compute $c^{(1)}_2$ we need $c_1^{(h)}$ for $h=1,\dots,p$. By \eqref{cnp} one can prove by induction that for $h\in \mathbb N$
$$c_1^{(h)}=-h u_0^h({a_0}-{a_1}u_0^p).$$
Hence we apply \eqref{c1n}  and we get
\begin{align*}
c^{(1)}_{2}&=-c_1^{(1)} \big({a_0} - {a_1}c^{(p)}_0 \big)+ {a_1}\left[\!\begin{matrix} \,1 \, \\ \, 1\,  \end{matrix}\!\right]_\beta \,\, c^{(p)}_1c^{(1)}_{0}
\end{align*}
from which we deduce by a direct computation
\begin{equation}\label{c12}c^{(1)}_{2}=u_0({a_0}-{a_1}u_0^p)({a_0}-{a_1}(p+1)u_0^p).\end{equation}

To compute  $c^{(1)}_3$ we also need $c_2^{(h)}$ for $h=1,\dots,p$.  Using \eqref{cnp} we obtain 
$$c^{(h)}_{2}=\sum_{k=0}^2 \left[\!\begin{matrix} \,2 \, \\ \, k\,  \end{matrix}\!\right]_\beta\, c^{(h-1)}_k\, \, c^{(1)}_{2-k}$$
and, by an inductive argument, using also the fact that $c_0^{(p)}=u_0^p$ for all $p\geq 1$, one can prove the closed formula for all $h\geq 1$
$$c^{(h)}_{2}= u_0^h({a_0}-{a_1}u_0^p)\left(h({a_0}-{a_1}(p+1)u_0^p)+\frac{h(h-1)}{2} \left[\!\begin{matrix} \,2 \, \\ \, 1\,  \end{matrix}\!\right]_\beta\, ({a_0}-{a_1}u_0^p)\right).$$
Then \begin{equation}\label{c13}
\begin{split}c^{(1)}_3=&-{a_0}c_2^{(1)}+{a_1}c_2^{(p+1)}\\
=&-{a_0}u_0({a_0}-{a_1}u_0^p)({a_0}-{a_1}(p+1)u_0^p)\\
&+{a_1}u_0^{p+1}({a_0}-{a_1}u_0^p)\bigg((p+1)({a_0}-{a_1}(p+1)u_0^p)\\
&\qquad\qquad\qquad\qquad\qquad+\frac{p(p+1)}{2} \left[\!\begin{matrix} \,2 \, \\ \, 1\,  \end{matrix}\!\right]_\beta\, ({a_0}-\
{a_1}u_0^p)\bigg)\\
=&-u_0({a_0}-{a_1}u_0^p)({a_0}-{a_1}(p+1)u_0^p)^2\\&+\frac{p(p+1)}{2} \left[\!\begin{matrix} \,2 \, \\ \, 1\,  \end{matrix}\!\right]_\beta\, {a_1}u_0^{p+1}({a_0}-{a_1}u_0^p)^2\\
=&-u_0({a_0}-{a_1}u_0^p)\bigg(({a_0}-{a_1}(p+1)u_0^p)^2\\
&\qquad\qquad\qquad-\frac{p(p+1)}{2} 
\left[\!\begin{matrix} \,2 \, \\ \, 1\,  \end{matrix}\!\right]_\beta\, {a_1}u_0^{p}({a_0}-{a_1}u_0^p)\bigg).
\end{split}\end{equation}

%
%
}

\section{Fractional logistic equations}

Here we extend  some of the results established in \cite{DL} for the fractional logistic equation with initial datum $u_0=1/2$ to the case of general initial data $u_0\in(0,1)$.
Applying the above method, if $p=1$ and $a_0=a_1=-1$ then the solution of 
\begin{align}
D^\beta_t u =  u -  u^{2}; \quad u(0)=u_0
\label{euleq}
\end{align}
can be represented in series form 
\begin{align*}
 u(t) = \sum_{n \geq 0} c^{(1)}_n \frac{t^{\beta n}}{\Gamma(\beta n +1)}
\end{align*}
where 
$$c_0^{(1)}=u_0; \qquad c^{(1)}_{n+1}=c^{(1)}_n-\sum_{k=0}^n\left[\!\begin{matrix} \, n \, \\ \, k\,  \end{matrix}\!\right]_\beta\, c^{(1)}_k\, \, c^{(1)}_{n-k}.$$
Note that, as shown in \cite{DL}, when $u_0=1/2$ above formula reduces to
$$\begin{cases}c_0^{(1)}=\dfrac{1}{2}; \quad c_1^{(1)}=\dfrac{1}{4};\\\\
  c^{(1)}_{2n}=0; \quad c^{(1)}_{2n+1}=-\displaystyle{\sum_{k=0}^n}\left[\!\begin{matrix} \, 2n \, \\ \, 2k+1\,  \end{matrix}\!\right]_\beta\, c^{(1)}_{2k+1}\, \, c^{(1)}_{2(n-k)-1}& \text{for }n\geq 1.\end{cases}$$

Keeping $a_0=a_1=-1$ and considering the equation with generic  $p$ we have that the solution of 
\begin{align}
D^\beta_t u =  u -  u^{p+1}; \quad u(0)=u_0
\label{euleq}
\end{align}
can be represented in series form 
\begin{align*}
 u(t) = \sum_{n \geq 0} c_n^{(p)} \frac{t^{\beta n}}{\Gamma(\beta n +1)}
\end{align*}
where
$$\begin{cases}
c^{(1)}_0 = u_0\quad c^{(1)}_1 =  u_0(1-u_0^p);\label{c11l}\\
\notag\\
 c^{(1)}_{n+1} = c_n^{(1)}- \displaystyle{\sum_{k=0}^n }\left[\!\begin{matrix} \, n \, \\ \, k\,  \end{matrix}\!\right]_\beta\, c^{(p)}_k\, \, c^{(1)}_{n-k}.
\end{cases}$$
\begin{remark}[On the null coefficients in the general case]Note that if $u_0=\sqrt[p]{1/(p+1)}$ then, in view of \eqref{c12},  $c_{2}^{(1)}=0$, in agreement with the case $p=1$. However it is not possible  to deduce that $c_{2n}^{(1)}\equiv 0$ as in the case $p=1$, because even for $p=2$, choosing $u_0=\sqrt{1/3}$ one can numerically verify that  $c_2^{(1)}=0$ and   $c_{4}^{(1)}\not=0$. Then one may look for some other generalization, for instance  imposing $c_{p+1}^{(1)}=0$ and guess whether $c_{n(p+1)}^{(1)}=0$ for some $n\geq 1$. Also in this case the answer is negative. Indeed, using the last equality \eqref{c13} to solve the equation  $c_3^{(1)}=0$ with respect to the initial datum $u_0$, we get by a direct computation that  if $u_0=\frac{12+6 \pi +\sqrt{144+96 \pi }}{24+18 \pi }$ and if $p=2$ then $c_3^{(1)}=0$, but  symbolic numerical computations yield $c_n^{(1)}\not=0$ for $n=3,\dots,30$.
 \end{remark}

Consider now the case $a_0=a_1=1$. Then the solution of 
\begin{align}
D^\beta_t u =  u -  u^{2}; \quad u(0)=u_0
\label{euleq}
\end{align}
can be represented as 
\begin{align*}
 u(t) = \sum_{n \geq 0} \bar c_n^{(p)} \frac{t^{\beta n}}{\Gamma(\beta n +1)}.
\end{align*}
where
\begin{align}
c^{(1)}_0 &= u_0\qquad c^{(1)}_1 =  -u_0(1-u_0^p);\label{c11l}\\
\notag\\
 c^{(1)}_{n+1} &= -c_n^{(1)}+ \sum_{k=0}^n \left[\!\begin{matrix} \, n \, \\ \, k\,  \end{matrix}\!\right]_\beta\, c^{(p)}_k\, \, c^{(1)}_{n-k}.
\end{align}

We compare the coefficients $c^{(1)}_n$ and $\bar c^{(1)}_n$. We have $\bar c^{(1)}_0 = c^{(1)}_0$, $\bar c^{(1)}_1 = -c^{(1)}_1$ and, by induction, $\bar c^{(p)}_1=-c^{(p)}_1$. Moreover $\bar c^{(1)}_2 = -c^{(1)}_2$.
However this symmetry breaks as soon as we consider $\bar c^{(2)}_2 $: indeed we have 
$$\bar c^{(2)}_2= - c^{(2)}_2+2 \left[\!\begin{matrix} \, 2 \, \\ \, 1\,  \end{matrix}\!\right]_\beta\,(c_1^{(1)})^2.$$

\subsection{Numerical simulations}
In our tests, we focused on the logistic case $a_0=a_1=-1$ and on the case $a_0=a_1=1$.  We computed the coefficients $c_n{(1)}$ using the recursive formulas \eqref{c11} and \eqref{cnp} and we approximated the solution of \eqref{eq1} with the partial sum 
$$u(t)=\sum_{n=0}^N c_n^{(1)}\frac{t^{\beta n}}{\Gamma(\beta n +1)}$$
with $N=200$ -- the parameter $N$ was tuned so that no appreciable difference can be noted with respect to higher order approximations. The method was validated by a comparison with the exact solutions of \eqref{eq1}, that can be explicitly computed in the ordinary case $\beta=1$.

\begin{figure}[h]
{\centering \includegraphics[scale=.45]{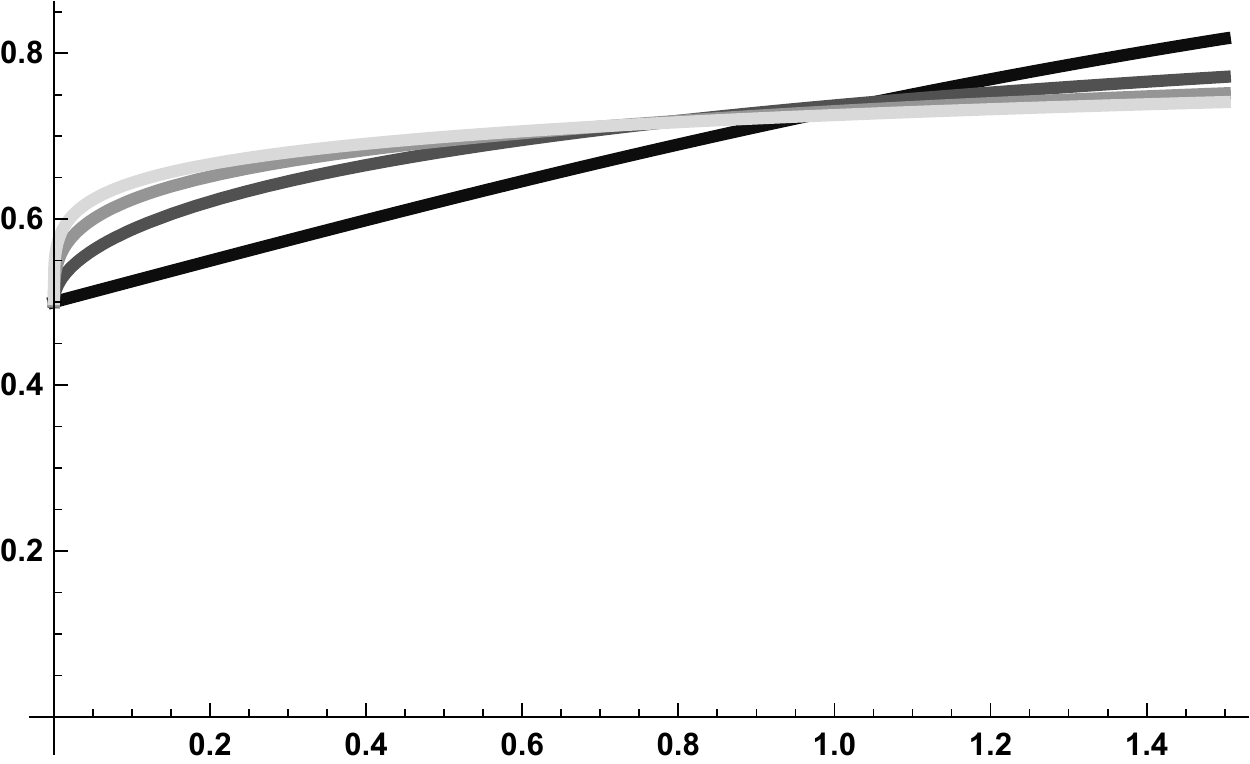}\quad\includegraphics[scale=.45]{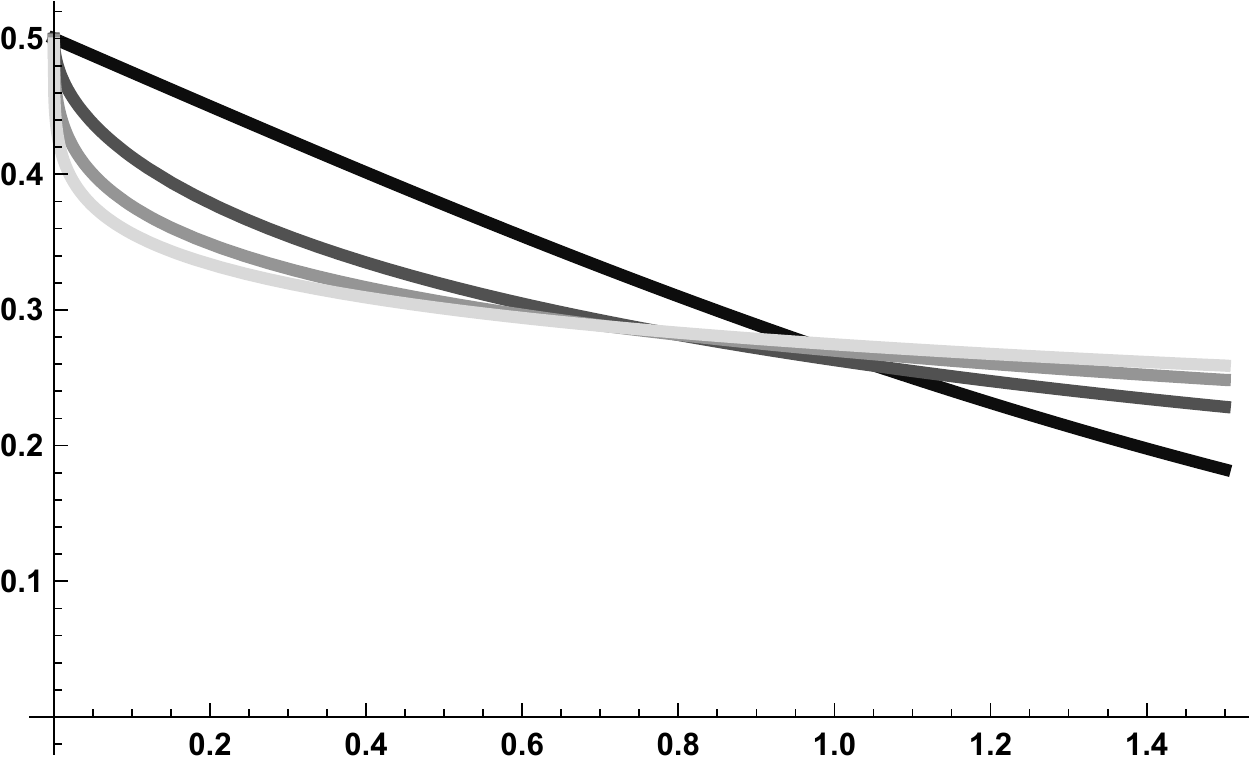}}
\caption{Numerical solutions of $D^\beta u-u=-u^2$ (on the left) and of $D^\beta u=-u+u^2$ (on the right) with $u_0=1/2$ and $\beta=1,1/2,1/3,1/4$. Darker lines correspond to bigger $\beta$'s. t\label{f2}}
\end{figure}

In Figure \ref{f2} we evaluated $u(t)$ with fixed initial datum $u_0=1/2$ and $p=1$ and several orders of fractional derivative: we may note the expected damping effect of fractional derivation.
\begin{figure}[h]
\centering \includegraphics[scale=.45]{logistic_datam.pdf}\quad  \includegraphics[scale=.45]{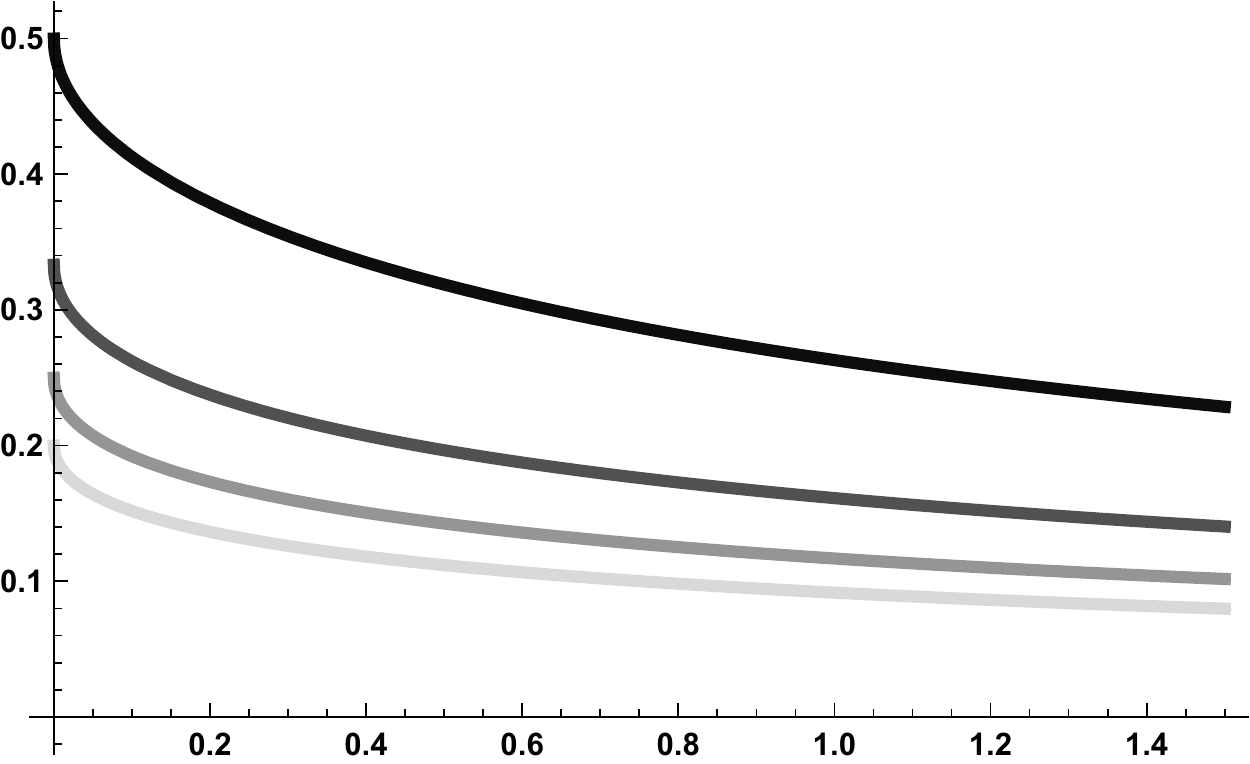}
\caption{Solutions of $D^\beta u-u=-u^2$ (on the left) and of $D^\beta u +u=u^2$ (on the right) with $\beta=1/2$ and $u_0=1/2,1/3,1/4,1/5$. Darker lines correspond to bigger $u_0$'s.\label{f3}}
\end{figure}

Figure \ref{f3} compares the solutions $u(t)$ with different initial data, setting $\beta=1/2$ and $p=1$. The resulting set of ordered curves suggest  local uniqueness of the solutions, whose investigation is however beyond the purpose of the present paper. 

\begin{figure}[h]
\centering \includegraphics[scale=.45]{euler_pm.pdf}\quad\includegraphics[scale=.45]{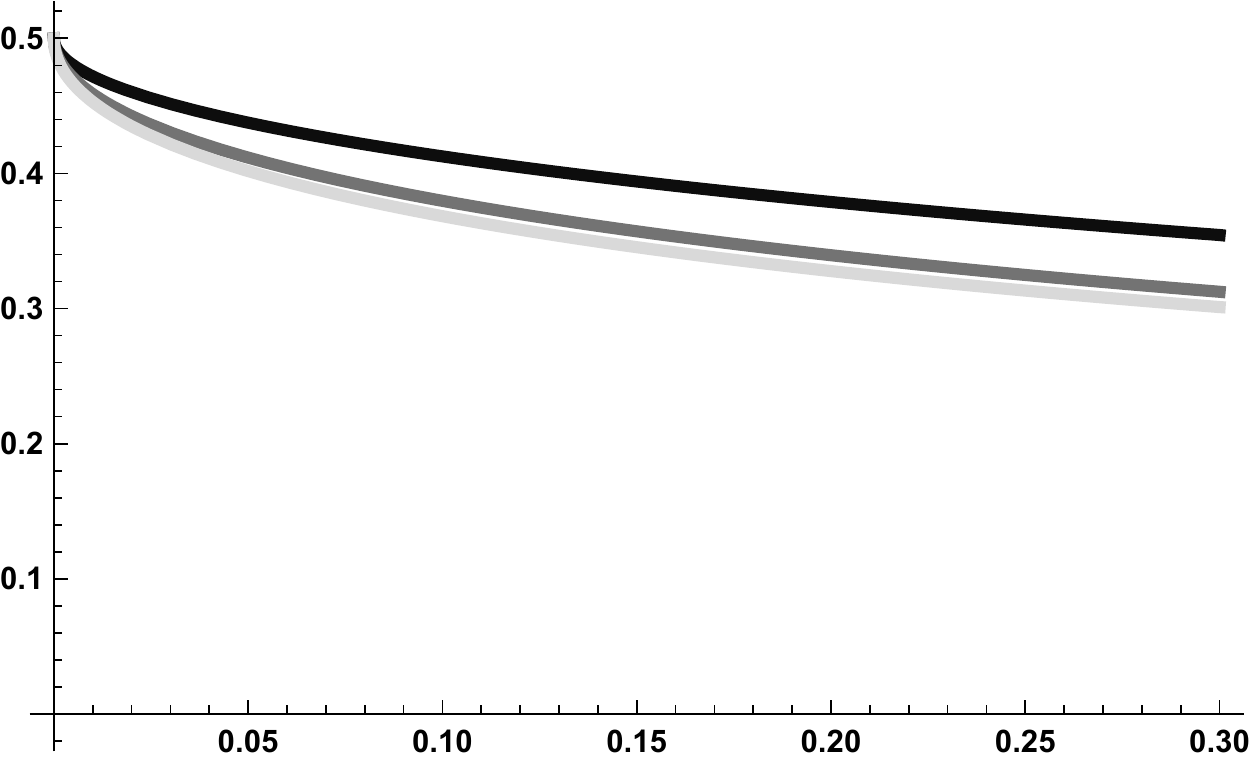}
\caption{Solutions of $D^\beta u-u=-u^{p+1}$ (on the left) and of $D^\beta u+u=u^{p+1}$ with $\beta=u_0=1/2$ and $p=1,2,3$. Darker lines correspond to bigger $p$'s.\label{f4}}
\end{figure}

We then investigated higher degree fractional Euler equations, setting $\beta=u_0=1/2$ and letting $p$ vary between 1 and 3, see Figure \ref{f4}. At least near $0$, from a qualitative point of view  the solutions display a similar behavior and no intersections between the solutions are detected. 

\begin{figure}[h]
\centering \includegraphics[scale=.45]{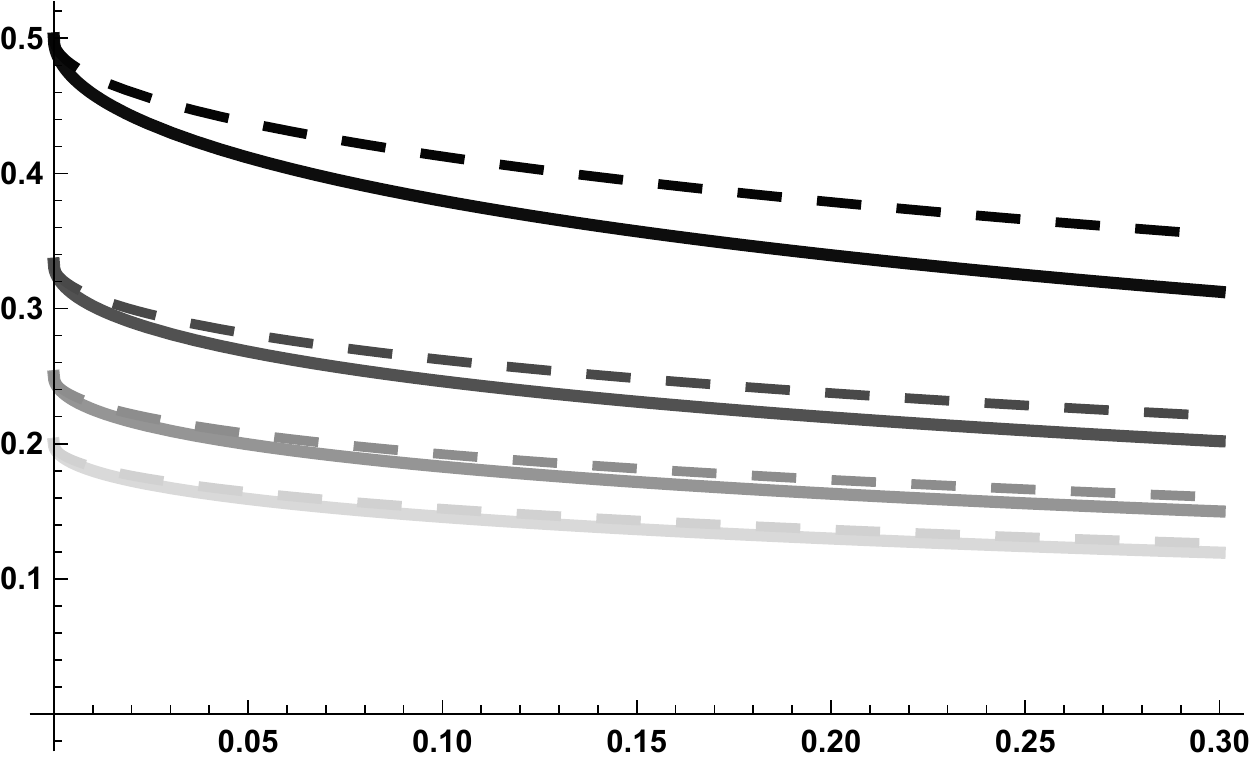}	\quad \includegraphics[scale=.45]{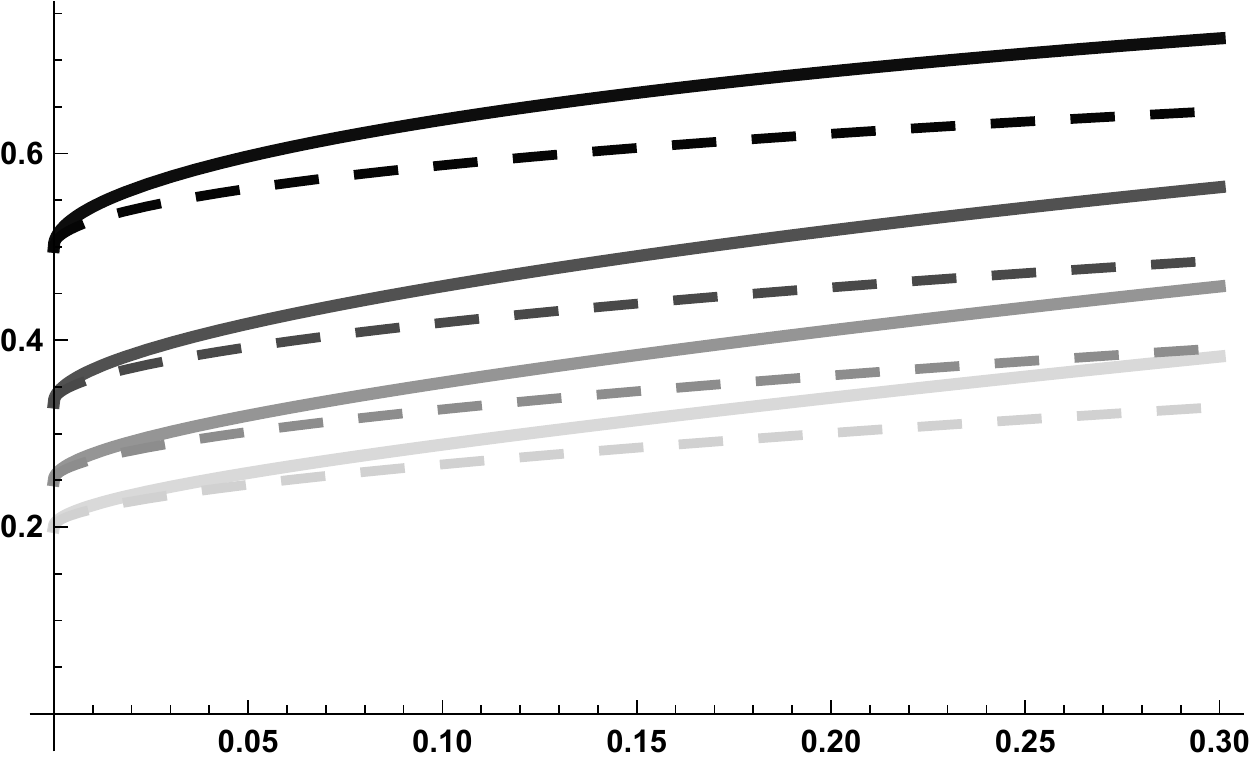}
\caption{Solutions of $D^\beta u+u=u^{p+1}$  (on the left) and of $D^\beta u=u-u^{(p+1)}$ (on the right)  with $\beta=1/2$ and $u_0=1/2,1/3,1/4,1/5$ and $p=1,2$. Darker lines correspond to bigger $u_0$'s; continuous lines correspond to $p=2$ and dashed lines correspond to $p=1$.\label{f5}}
\end{figure}

For the seek of comparison, we collected some of the above results in Figure \ref{f5}, showing the combined effect of varying initial data $u_0$ and degrees $p$. 

\begin{figure}[h]
\centering \includegraphics[scale=.45]{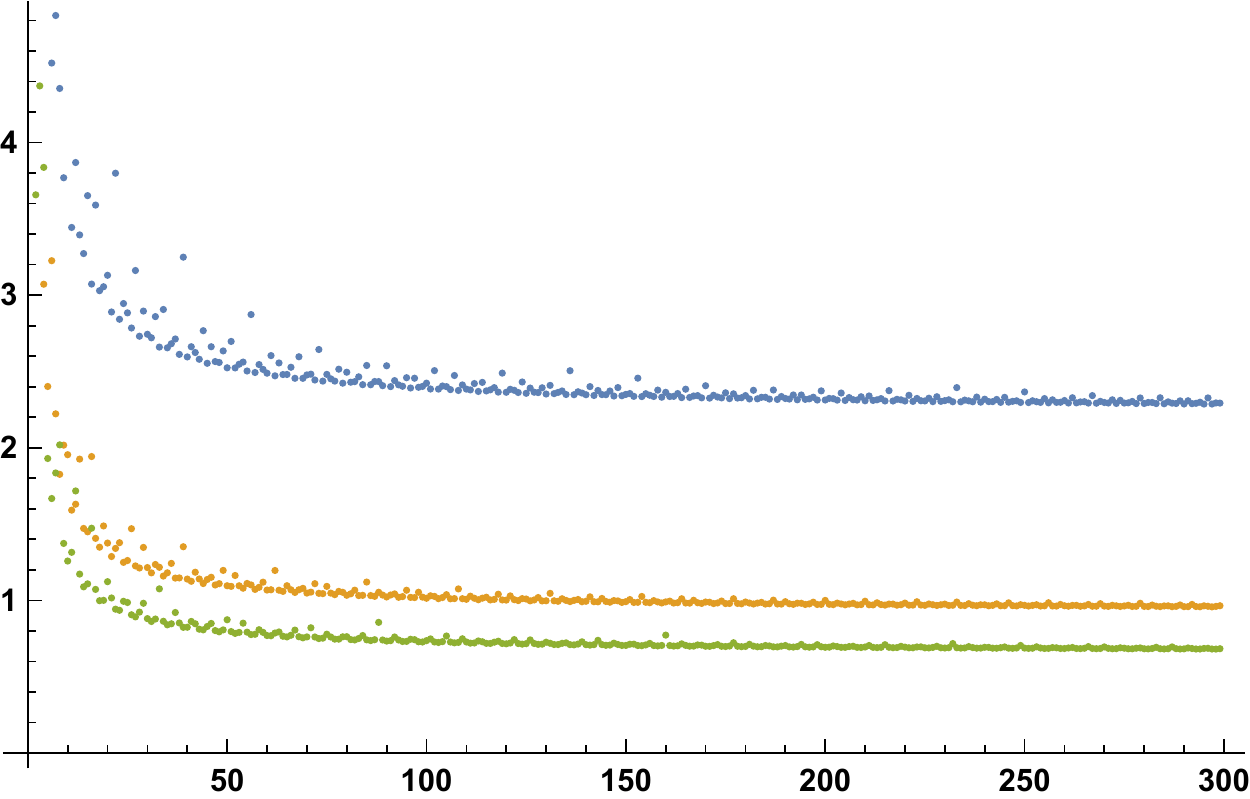}\quad \includegraphics[scale=.45]{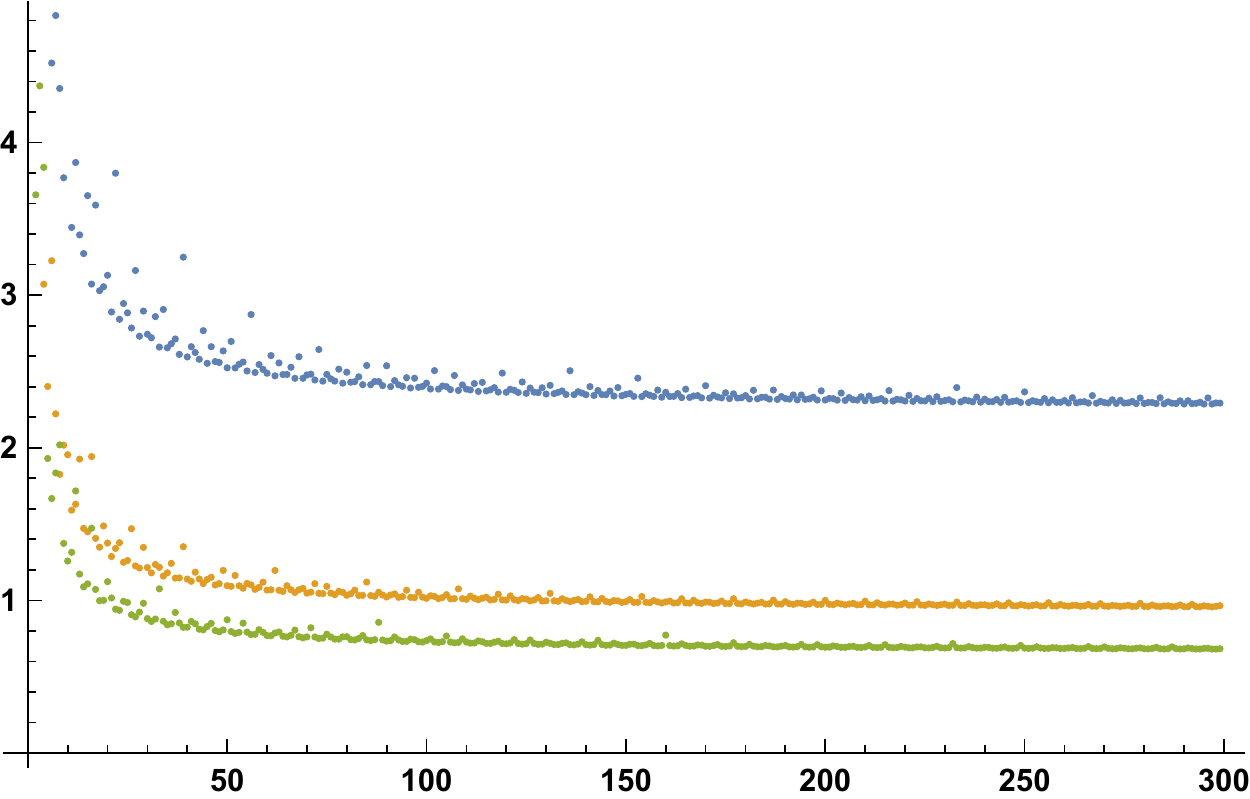}
\caption{Asymptotic behavior of the sequence approximating the radius of convergence $r_n$ for $a_0=a_1=-1$ (on the left) and $a_0=a_1=1$ (on the right)  with $\beta=1/2$, $u_0=1/3$ and $p=1,2,3$ in blue, orange, green respectively.\label{f6}}
%
\end{figure}

\begin{figure}[h]
\centering \includegraphics[scale=.45]{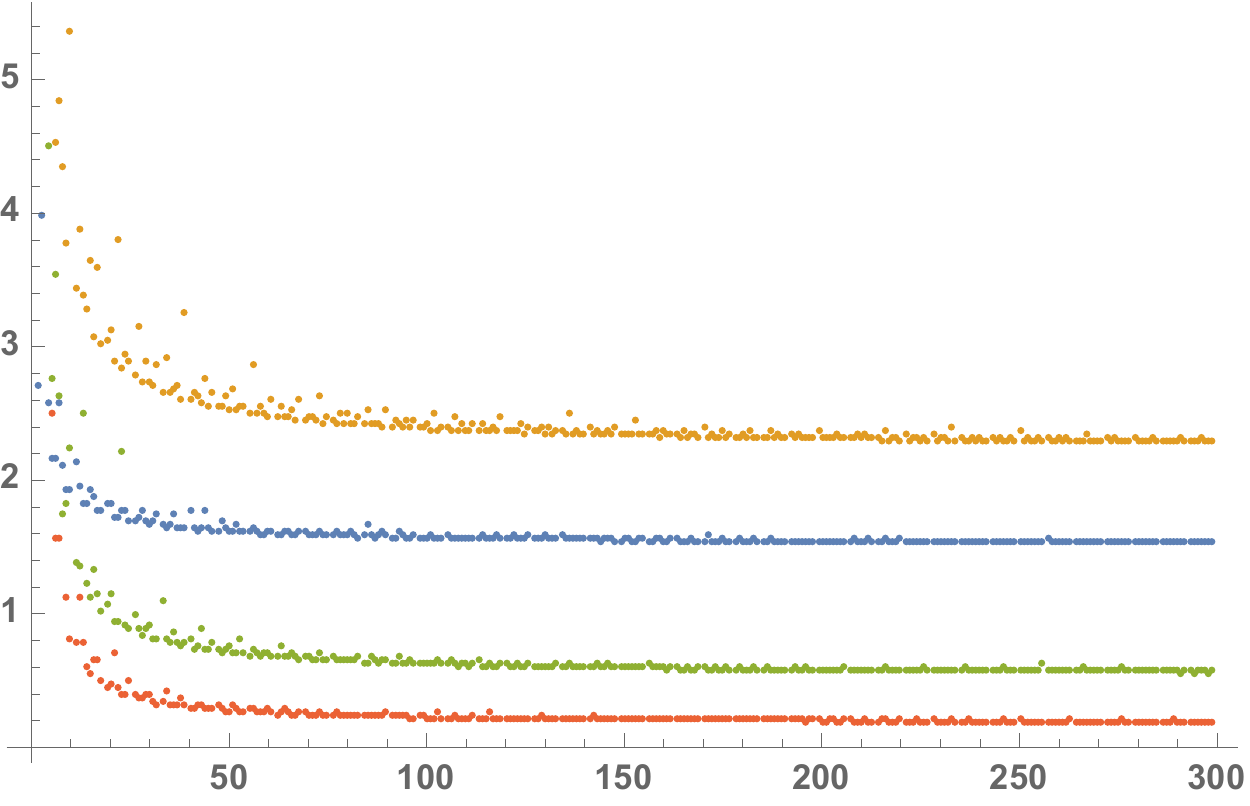}\quad \includegraphics[scale=.45]{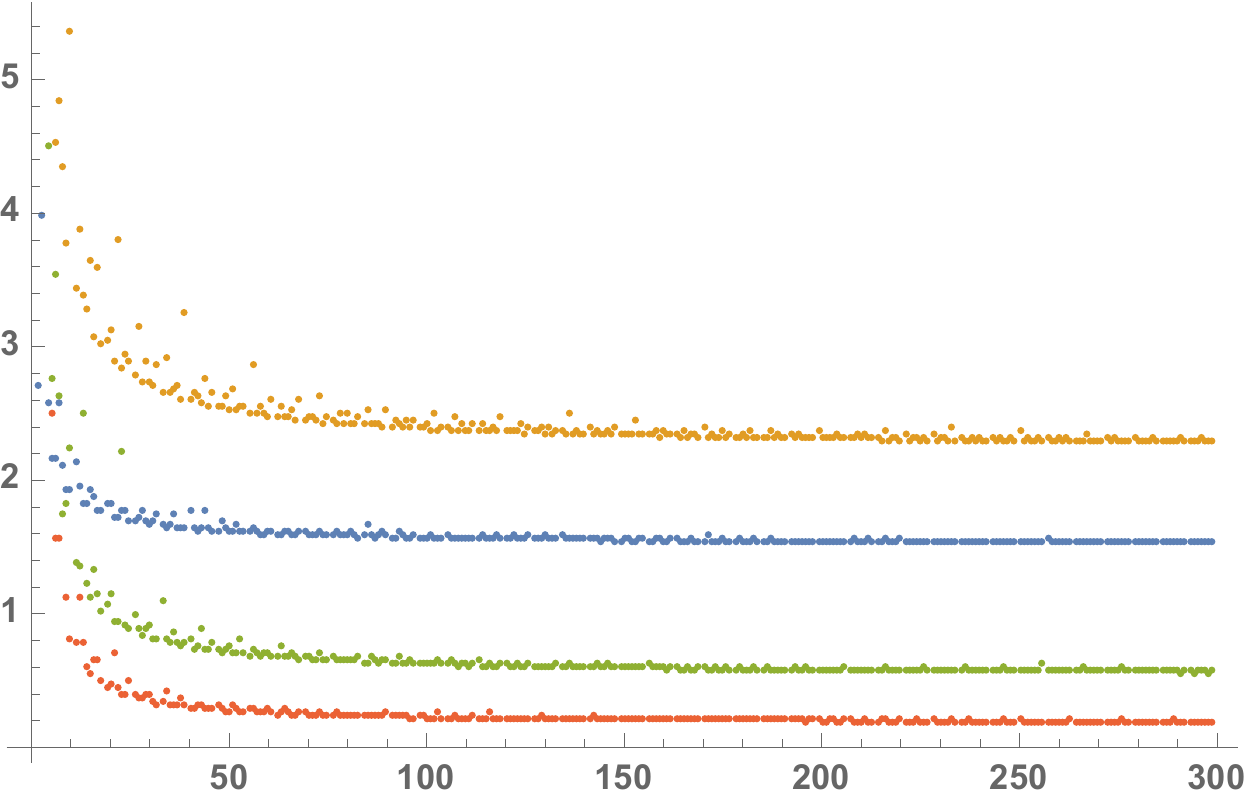}
\caption{Asymptotic behavior of the sequence approximating the radius of convergence $r_n$ for $a_0=a_1=-1$ (on the left) and $a_0=a_1=1$ (on the right)  with $u_0=1/3$, $p=1$, and $\beta=1,1/2,1/3,1/4$ in blue, orange, green and red respectively.\label{f7}}
\end{figure}

Finally we propose some numerical estimations for the radius of convergence of the series \eqref{series}, by computing the sequence
$$r_n:=\left(\frac{\Gamma(\beta n+1)}{|c_n^{(1)}|}\right)^{\frac{1}{\beta n}}.$$
In Figure \ref{f6} and Figure \ref{f7} we plotted the first $300$ terms of $r_n$ with varying degrees $p$ and orders of derivation $\beta$. The asymptotic behavior of $r_n$ suggests an exponential increase for the series coefficients $c_n^{(1)}$ in all the cases under exam. Furthermore, their comparison shows  the radius of convergence $r:=\lim r_n$ to 
be decreasing with respect to both the degree $p$ (Figure \ref{f6}) and order of derivation $\beta$ (Figure \ref{f7}). Finally, no substantial difference betweens the case $a_0=a_1=-1$ and the case $a_0=a_1=1$ emerged. 

\end{document}